\documentclass[a4paper]{article}
\usepackage{amsmath,amsfonts,graphicx}
\usepackage[all]{xy}
\newtheorem{prop}{Proposition}[section]
\newtheorem{defi}[prop]{Definition}
\newtheorem{lemma}[prop]{Lemma}

\newtheorem{teo}[prop]{Theorem}

\begin{document}

\begin{center}
\large{\bf {THE ALGEBRAS OF SEMI-INVARIANTS OF EUCLIDEAN QUIVERS}}\\

\bigskip

\large{Cristina Di Trapano}\\

\bigskip

\small{University of Rome "Tor Vergata"\\
Department of Mathematics\\
Via della Ricerca Scientifica 1, 00133 Rome (Italy)}\\

\bigskip

\verb"e-mail": cristinaditrapano@virgilio.it\\
\end{center}

\bigskip

\noindent ABSTRACT. We give a new short proof of Skowro\'nski and
Weyman's theorem about the structure of the algebras of
semi-invariants of Euclidean quivers, in the case of quivers
without oriented cycles and in characteristic zero. Our proof is
based essentially on Derksen and Weyman's result about the
generators of these algebras and pro\-per\-ties of Schofield
semi-invariants.\\

\bigskip

\bigskip
\begin{center}
\noindent{\Large{\textbf{{Introduction}}}}\\
\end{center}
\bigskip

Let $Q$ be a connected and finite quiver. Define
$Rep(Q,\underline{\alpha})$ the affine variety of representations
of $Q$ of dimension vector $\underline{\alpha}$. We are interested
in the action of the group $SL(\underline{\alpha}):= \prod_{x\in
Q_0}{SL(\underline{\alpha}(x))}$ on this variety. In particular,
we look at the algebra of semi-invariants
$SI(Q,\underline{\alpha}):=K[Rep\,(Q,\underline{\alpha})]^{SL(\underline{\alpha})}$.
By Sato-Kimura's lemma ([SK], Section 4 Proposition 5) it follows that for Dynkin qui\-vers the algebra of semi-invariants is a polynomial algebra.\\
In ([SW], Theorem 1), Skowro\'nski and Weyman prove the following
theorem:
\begin{teo}
For each dimension vector $\underline{\alpha}$, the algebra
$SI(Q,\underline{\alpha})$ is either a polynomial algebra or a
quotient of a polynomial algebra by a principal ideal if and only
if $Q$ is a Dynkin quiver or a Euclidean quiver.
\end{teo}

\noindent In particular, they give an explicit description of
$SI(Q,\underline{d})$ where $Q$ is a Euclidean quiver (also with
oriented cycles) and $\underline{d}$ a dimension vector ([SW],
Theo\-rem 21). On the other hand in ([Sh], Theorem 8.6) Shmelkin
gives an independent proof of Theorem 0.1 based on some Luna's
results ([Lu1], [Lu2]).

In this paper, assuming that $Q$ is a quiver without oriented
cycles and $K$ is an algebraically closed field of characteristic
zero, we find the same presentation of the algebra
$SI(Q,\underline{d})$, $Q$ a Euclidean quiver and $\underline{d}$
a dimension vector, in the sense described in [SW], using new
short methods. Our proof is based on Derksen and Weyman's theorem
about generators of algebras of semi-invariants for quivers
without oriented cycles ([DW], Theorem 1), on some properties of
Schofield semi-invariants ([DW], Lemma 1) and on Derksen,
Schofield and Weyman's theorem relating the dimension of
weight spaces of semi-invariants to the number of subrepresentations of a general quiver representation ([DSW], Theorem 1).\\
The paper is organized as follows. In Section 1, we recall some
results about Schofield semi-invariants. In Section 2, we recall
some facts about Euclidean quivers. In Section 3, we formulate the
results on the structure of the algebras of semi-invariants of
Euclidean quivers. In Section 4, we provide the proofs of the
results stated in Section 3.\\

{\bf{Acknowledgments:}} I would like to thank F. Gavarini, E.
Strickland and J. Weyman for helpful discussions.\\

\section{Preliminary results}

Let $K$ be an algebraically closed field of characteristic zero.\\
A \textbf{quiver} is a directed graph $Q = (Q_0, Q_1)$, where
$Q_0$ is the set of vertices and $Q_1$ is the set of arrows. Each
arrow has its head $ha$ and tail $ta$, both in $Q_0$;
$$ta\stackrel{a}{\longrightarrow}ha.$$

\medskip

\noindent We assume that the quiver $Q$ is ${\bf{finite}}$ (that
is, $Q_0$ and $Q_1$ are finite) and ${\bf{connected}}$.\\
A ${\bf{path}}$ \emph{p} is a sequence $a_1a_2\cdots a_m$ of
arrows such that $ta_i = ha_{i+1}$ for $i = 1,\ldots, m-1$. For
every vertex $x\in Q_0$ we also have a trivial path $e_x$ with
$he_x = te_x = x$. We say that $Q$ has no oriented cycles if there
is no nontrivial path $p$ such that
$tp = hp$.\\
We will assume from now on that $Q$ has $\bf{no}$ $\bf{oriented}$ $\bf{cycles}$.\\

A \textbf{representation} $V$ of $Q$ is a collection
$(V(x)\;|\;x\in Q_0)$ of finite dimensional $K$-vector spaces
together with a collection of $K$-linear maps $(V(a):
V(ta)\rightarrow V(ha)\;|\;a\in Q_1)$. For every representation
$V$ we define the \textbf{dimension vector} $\underline{d}_V :
Q_0\rightarrow \mathbb{N}$ of $V$ by
$\underline{d}_V(x) := dim_K V(x),\; x\in Q_0.$\\
Suppose that $V$ and $W$ are representations of a quiver $Q$. A
\textbf{morphism} $f :V\rightarrow W$ is a collection of
$K$-linear maps $(f(x): V(x)\rightarrow W(x)\;|\; x\in Q_0)$ such
that for each arrow $a\in Q_1$ $$f(ha)V(a)=W(a)f(ta).$$ If
moreover, $f(x)$ is invertible for each $x\in Q_0$, then $f$ is
called an \textbf{isomorphism}. We denote the linear space of
morphisms from $V$ to $W$ by $Hom_Q
(V, W)$.\\

Let $Q$ be a quiver as above and $\underline{\alpha}$ a dimension
vector. We can view a representation of $Q$ of dimension vector
$\underline{\alpha}$ as an element in
$$Rep\,(Q, \underline{\alpha}) := \bigoplus_{a\in Q_1}{Hom(K^{\underline{\alpha}(ta)},
K^{\underline{\alpha}(ha)})}.$$ The group $GL(\underline{\alpha})
:= \prod_{x\in Q_0}{GL(\underline{\alpha}(x))}$ and its subgroup
$SL(\underline{\alpha})$ act on the re\-pre\-sen\-ta\-tion space
$Rep\;(Q,
\underline{\alpha})$ in the following way:\\

$$(\phi\cdot V)(a) :=\phi(ha)V(a)\phi(ta)^{-1}$$\\

\noindent where $\phi = (\phi(x) \in
GL(\underline{\alpha}(x))\;|\; x\in Q_0) \in
GL(\underline{\alpha})$ and $V \in Rep(Q, \underline{\alpha})$. We
will look at the ring of $SL(Q,\underline{\alpha})$-invariants
which is isomorphic to \textbf{the ring of semi-invariants}
$$SI(Q, \underline{\alpha}) := \bigoplus_{\chi\in char (GL(Q,\underline{\alpha}))}{SI(Q,
\underline{\alpha})_\chi}$$ where
$$SI(Q,\underline{\alpha})_\chi :=\{f\in K[Rep\,(Q, \underline{\alpha})]\; |\; g\cdot f = \chi(g)f ,\; g \in
GL(\underline{\alpha})\}.$$ Suppose that $\chi :
GL(\underline{\alpha})\rightarrow K^*$ is a character. Such
character always looks like
$$(\phi(x)\in GL(\underline{\alpha}(x))\;|\;x\in Q_0)\rightarrow \prod_{x\in
Q_0}{det(\phi(x))^{\sigma(x)}}.$$ Here any map
$\sigma:Q_0\rightarrow
\mathbb{Z}$ is called a \textbf{weight}.\\

Next, we recall Schofield semi-invariants [Sc], as they are
the main objects we shall use in our proofs.\\
For representations $V$ and $W$ of $Q$ there is a canonical exact
sequence ([Ri1]):

\begin{align}
\nonumber 0\rightarrow Hom_Q(V,W)\rightarrow \bigoplus_{x\in
Q_0}{Hom(V(x),W(x))}&\buildrel{d^W_V}\over\rightarrow\\
\nonumber\buildrel{d^W_V}\over\rightarrow\bigoplus_{a\in
Q_1}{Hom(V(ta),W(ha))}&\rightarrow Ext_Q(V,W)\rightarrow 0.
\end{align}

\noindent The map $d_V^W$ is defined by
$$(\phi(x)\;|\;x\in Q_0)\mapsto (W(a)\phi(ta)-\phi(ha)V(a)\;|\;a\in
Q_1).$$ For $\underline{\alpha}, \underline{\beta}$ dimension
vectors, we define the $\textbf{Euler form}$
$$\langle{\underline{\alpha},\underline{\beta}}\rangle =
\sum_{x\in
Q_0}{\underline{\alpha}(x)\underline{\beta}(x)}-\sum_{a\in
Q_1}{\underline{\alpha}(ta)\underline{\beta}(ha)}.$$ Let us choose
the dimension vectors $\underline{\alpha}$ and $\underline{\beta}$
such that $\langle{\underline{\alpha},\underline{\beta}}\rangle
=0$. Then for every $V\in Rep(Q,\underline{\alpha})$ and $W \in
Rep(Q,\underline{\beta})$ the matrix of $d_W^V$ is a square
matrix. Following [Sc], we can define the semi-invariant $c \in
K[Rep(Q, \underline{\alpha})\times Rep(Q, \underline{\beta})]$ by
$c(V, W):=
det(d^V_W)$.\\
\noindent For a fixed $V$, the restriction of $c$ to $\{V\}\times
Rep(Q,\underline{\beta})$ defines a semi-invariant
$$c^V := c(V,-) \in
K[Rep(Q,\underline{\beta})]^{SL(\underline{\beta})}=SI(Q,\underline{\beta}).$$
\noindent Similarly, for a fixed $W$, the restriction of $c$ to $
Rep(Q,\underline{\alpha})\times \{W\}$ defines a semi-invariant
$$c_W := c(-,W) \in
K[Rep(Q,\underline{\alpha})]^{SL(\underline{\alpha})}=SI(Q,\underline{\alpha}).$$\\
\noindent The semi-invariants
$c^V$ and $c_W$ are called \textbf{Schofield semi-invariants} corresponding to $V$ and $W$ respectively.\\
These semi-invariants have the following important properties ([Sc] Lemma 1.4; [DW] Lemma1):\\

\begin{prop}
\item{a)} The semi-invariant $c^V$ lies in
$SI(Q,\underline{\beta})_{\langle{\underline{\alpha},-}\rangle}$
for every $V\in
Rep(Q,\underline{\alpha})$.\\

\item{b)} The semi-invariant $c_W$ lies in
$SI(Q,\underline{\alpha})_{-\langle{-,\underline{\beta}}\rangle}$
for every $W\in
Rep(Q,\underline{\beta})$.\\

\item{c)} Let
$$0\rightarrow V'\rightarrow V\rightarrow V''\rightarrow 0$$
be an exact sequence of representations, with
$\underline{\alpha}':= \underline{d}_{V'}$. If
$\langle{\underline{\alpha}',\underline{\beta}}\rangle = 0$, then
$c^V=c^{V'}c^{V''}$. If
$\langle{\underline{\alpha}',\underline{\beta}}\rangle < 0$,
then $c^V = 0$.\\

\item{d)} Let
$$0\rightarrow W'\rightarrow W\rightarrow W''\rightarrow 0$$
be an exact sequence of representations, with
$\underline{\beta}':= \underline{d}_{W'}$. If
$\langle{\underline{\alpha},\underline{\beta}'}\rangle = 0$, then
$c_W =c_{W'}c_{W''}$. If
$\langle{\underline{\alpha},\underline{\beta}'}\rangle > 0$,
then $c_W = 0$.\\

\item{e)}
$c^V(W)=0\Leftrightarrow\;Hom_Q(V,W)\neq0\Leftrightarrow\;Ext_Q(V,W)\neq0.$\quad$\diamondsuit$
\end{prop}
{\ }\\

\noindent Let's recall two important results on semi-invariants
which our proofs are based on.

\begin{teo}\emph{([DW] Theorem 1, [Ch] Corollary 2.5)} Let $Q$ be a quiver without
oriented cycles, $\underline{\beta}$ a dimension vector, $\sigma$
a weight.
\item{a)} For any dimension vector $\underline{\gamma}$, a weight space
$SI(Q,\underline{\gamma})_\sigma$ can be nonzero only for weights
satisfying $\sigma(\underline{\gamma})=\sum_{x\in
Q_0}{\sigma(x)\underline{\gamma}(x)}=0$.
\item{b)} If there is no dimension vector $\underline{\alpha}$ such that $\sigma$ and $\langle{\underline{\alpha},-}\rangle$ determine
the same character of $GL(\underline{\beta})$, then
$SI(Q,\underline{\beta})_{\sigma}=0$.
\item{c)} If $\sigma =\langle{\underline{\alpha},-}\rangle$ with
$\langle{\underline{\alpha},\underline{\beta}}\rangle=0$, then
$SI(Q,\underline{\beta})_\sigma$ is spanned as a vector space by
the semi-invariants $c^V$ for $V \in
Rep(Q,\underline{\alpha})$.\qquad$\diamondsuit$
\end{teo}

\noindent It follows from Proposition 1.1 and Theorem 1.2 that
$SI(Q,\underline{\beta})$ is generated as a $K$-algebra by
Schofield semi-invariants $c^V$, for $V$ indecomposable
representation of $Q$ with
$\langle{\underline{d}_V,\underline{\beta}}\rangle = 0$.\\
Before giving the next result, we recall the notion of general
representation ([Ka],[Ka1] Section 4).\\
We say that a \textbf{general representation} with dimension
vector $\underline{d}$ has a certain property, if all
representations in some Zariski open (and dense) subset of the
space of $\underline{d}$-dimensional representations have that
property. We say that
$$
\underline{d}=\underline{d}_1+\underline{d}_2+\ldots+\underline{d}_r
$$
is the generic decomposition of a dimension vector $\underline{d}$
if a general representation $W$ of dimension vector
$\underline{d}$ has a decomposition $W = W_1 \oplus
W_2\oplus\ldots\oplus W_r$, where
each $W_i$ is indecomposable of dimension vector $\underline{d}_i$.\\

\begin{teo}\emph{([DSW] Theorem 1)}
Let $Q$ be a quiver without oriented cycles and
$\underline{\alpha}$ and $\underline{\beta}$ two dimension
vectors. Let $N(\underline{\beta}, \underline{\alpha})$ be the
number of $\underline{\beta}$-dimensional subrepresentations of a
general $\underline{\alpha}$-dimensional representation, and
$M(\underline{\beta},\underline{\alpha})$ the dimension of the
space of semi-invariants of weight
$\langle{\underline{\beta},-}\rangle$ on the representation space
of dimension vector
$\underline{\gamma}:=\underline{\alpha}-\underline{\beta}$. If
$\langle{\underline{\beta},\underline{\gamma}}\rangle=0$, then
$N(\underline{\beta},
\underline{\alpha})=M(\underline{\beta},\underline{\alpha})$.\;$\diamondsuit$
\end{teo}
{\ }\\

\section{Euclidean quivers}

\noindent We recall some facts about Euclidean quivers.\\
Let $Q$ be a Euclidean quiver without oriented cycles of type
$\tilde{A}_n, \tilde{D}_n, \tilde{E}_6, \tilde{E}_7, \tilde{E}_8$,
i.e a quiver for which the underlying graph is one of
the following type:\\

$$\tilde{A}_n :\begin{array}{ccccccccc}
   &  & c_1 & \makebox[0.5cm][c]{---} & \cdots & \makebox[0.5cm][c]{---}  & c_u &  &  \\
   & / &  &  &  &  &  & \backslash &  \\
  & a &  &  &  &  &  & b &  \\
   & \backslash &  & &  &  &  & / &  \\
   &  & d_1 & \makebox[0.5cm][c]{---}  & \cdots & \makebox[0.5cm][c]{---} & d_v &  &  \\
\end{array}
$$

where $u+v=n-1$,\\

$$\tilde{D}_n :
\begin{array}{ccccccccc}
   a_1 &  &  &  &  &  &  &  & b_1  \\
   & \backslash &  &  &  &  &  & / &  \\
   &  & z_1 & \makebox[0.5cm][c]{---} & \ldots & \makebox[0.5cm][c]{---} & z_{n-3} &  &  \\
  & / &  &  &  &  &  & \backslash &  \\
  a_2 & &  &  &  &  &  &  & b_2 \\
\end{array}
$$

$$\tilde{E}_6 :
\begin{matrix}
   &  &  &  & c_1 &  &  &  & \cr
   &  &  &  & |& &  &  &  \cr
   &  &  &  & c_2 &  &  &  &  \\
   &  &  &  & |& &  &  &  \\
  a_1 & \makebox[0.5cm][c]{---} & a_2 & \makebox[0.5cm][c]{---} & z & \makebox[0.5cm][c]{---} & b_2 & \makebox[0.5cm][c]{---} & b_1 \\
\end{matrix}
$$

\bigskip

\bigskip

$\tilde{E}_7 :\begin{array}{ccccccccccccc}
   &  &  &  & c &  &  &  & & & & & \\
   &  & &  & | &  & & &  & & & & \\
  b_1 & \makebox[0.5cm][c]{---} & b_2 & \makebox[0.5cm][c]{---} & b_3 & \makebox[0.5cm][c]{---} & z &
  \makebox[0.5cm][c]{---} & a_3 & \makebox[0.5cm][c]{---} & a_2 &
  \makebox[0.5cm][c]{---}& a_1\\
\end{array}$

\bigskip

\bigskip

$\tilde{E}_8 :\begin{array}{ccccccccccccccc}
   &  &  &  & c &  &  &  & & & & & \\
   &  &  &  & | &  & & &  & & & & \\
  b_1 & \makebox[0.5cm][c]{---} & b_2 & \makebox[0.5cm][c]{---} & z & \makebox[0.5cm][c]{---} & a_5 & \makebox[0.5cm][c]{---} & a_4 &
  \makebox[0.5cm][c]{---}& a_3 &
  \makebox[0.5cm][c]{---}& a_2 &\makebox[0.5cm][c]{---}& a_1\\
\end{array}$\\

\bigskip
\bigskip

\noindent By [DR] Proposition 1.2, the quadratic form
$q_Q:\mathbb{Z}^{Q_0}\rightarrow\mathbb{Z}$ defined by
$$
q_Q(\underline{\alpha}):=\sum_{v\in
Q_0}{\underline{\alpha}(v)^2}-\sum_{a\in
Q_1}{\underline{\alpha}(ta)\underline{\alpha}(ha)}
$$
is positive semi-definite and there exists a unique dimension
vector $\underline{h} \in \mathbb{N}^{Q_0}$ such that
$\mathbb{Z}\underline{h}$ is the radical of $q_Q$ ([DR] page 9). \\
Following [DR] Section 1, we define the \textbf{defect} of a
module $V$ as $\partial(\underline{d}_V):=\langle{\underline{h},
\underline{d}_V}\rangle$, we say that an indecomposable
representation $V$ is \emph{preprojective, re\-gu\-lar} or
\emph{preinjective} if and only if the defect of $V$ is negative,
zero or positive, respectively. As we shall see later (Section
4.1), we are interested
in regular modules, so we give more details about them.\\

\noindent Regular representations of $Q$ form an abelian category
$Reg_K(Q)$ ([DR] Proposition 3.2), so we may speak of regular
composition series, simple regular objects, etc., referring to
composition series, simple objects, etc. inside the category
$Reg_K(Q)$. The category $Reg_K(Q)$ is serial ([DR] Theorem 3.5):
any indecomposable regular representation has a unique regular
composition series, thus it is
uniquely determined by its regular socle (simple regular) and by its regular length.\\
By [DR] Theorem 3.5, the category $Reg_K(Q)$ decomposes into a
direct sum of categories $R_t$, with $t\in K\cup
\{\infty\}=\mathbb{P}_1(K)$. We call each category $R_t$ a
$\textbf{tube}$. In order to describe such tubes, we need the
following definition:
\begin{defi}
Let $Q$ be a quiver without oriented cycles.  We may assume that
$Q_0 = \{1, 2,\ldots,n\}$ and for every $a\in Q_1$ we have
$ta<ha$. We define the \textbf{Coxeter element}
$$ C := \sigma_1\sigma_2\ldots\sigma_n$$
where each $\sigma_i$ acts on dimension vectors as follows:
$$\sigma_i(\underline{\alpha} ) (j)= \begin{cases}
                                            \underline{\alpha}(j)&if\;j\,\ne i;\\
                                            \sum_{a\in Q_1; ta=i} \underline{\alpha}(ha)+
                                            \sum_{a\in Q_1; ha=i}\underline{\alpha}(ta)-\underline{\alpha}(i),&otherwise.
                                            \end{cases}$$

\end{defi}
{\ }\\

\bigskip

\noindent By [DR] Lemma 1.3 and Lemma 3.3, for each simple regular
representation $V$, the orbit of the dimension vector of $V$ under
$C$ is always finite. In particular, a simple regular
representation $V$ with the dimension vector which is fixed by $C$
is called $\textbf{homogeneous}$. In such a case, we have that the
dimension vector of $V$ is equal to $\underline{h}$. About the
$C$-orbits of simple non homogeneous modules, we recall
the following description ([DR], Section 6):\\

\begin{prop}
Let $Q$ be a Euclidean quiver. Then there are at most three
$C$-orbits $\Delta = \{\underline{e}_i, i\in I\},\ \Delta'=
\{\underline{e}_i', i\in I'\},\ \Delta''= \{\underline{e}_i'',
i\in I''\}$, of dimension vectors of non homogeneous simple
regular representations of $Q$ ($I, I', I''$ could be empty). We
can assume $I = \{0,1,\ldots,u-1\},\ I'=\{0,1,\ldots,v-1\},\
I''=\{0,1,\ldots,w-1\}$ and
$C(\underline{e}_i)=\underline{e}_{i+1}$ for $i\in I\;
(\underline{e}_u=\underline{e}_0)$,\
$C(\underline{e}'_i)=\underline{e}'_{i+1}$ for $i\in I'\;
(\underline{e}'_v=\underline{e}'_0)$,\
$C(\underline{e}''_i)=\underline{e}''_{i+1}$ for $i\in I''\;
(\underline{e}''_w=\underline{e}''_0)$.\quad$\diamondsuit$
\end{prop}

\noindent Graphically we may represent them as the following polygons:\\

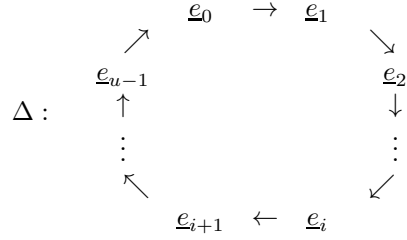
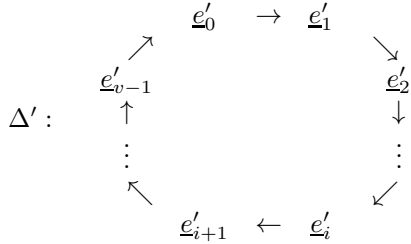
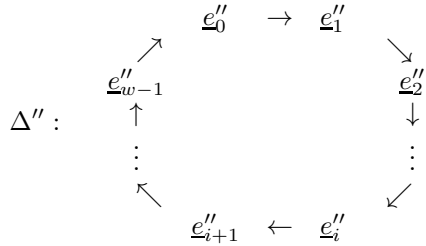
$$\Delta:
\begin{array}{ccccccc}
    &  & \underline{e}_0 & \rightarrow & \underline{e}_1 &  &  \\
   & \quad\nearrow &  &  & & \searrow &  \\
   & \underline{e}_{u-1} &  &  &  & \quad \underline{e}_2 &  \\
   & \uparrow &  &  &  & \quad\downarrow &  \\
   & \vdots &  &  &  & \quad\vdots &  \\
   & \quad \nwarrow&  &  &  & \swarrow &  \\
   &  & \underline{e}_{i+1} & \leftarrow & \underline{e}_i &  &  \\
\end{array}$$\\

$$\Delta':
\begin{array}{ccccccc}
    &  & \underline{e}_0' & \rightarrow & \underline{e}_1' &  &  \\
   & \quad\nearrow &  &  &  & \searrow &  \\
   & \underline{e}_{v-1}' &  &  &  & \quad \underline{e}_2' &  \\
   & \uparrow &  &  &  & \quad\downarrow &  \\
   & \vdots &  &  &  & \quad\vdots &  \\
   & \quad \nwarrow&  &  &  & \swarrow & \\
   &  & \underline{e}_{i+1}' & \leftarrow & \underline{e}_i' &  &  \\
\end{array}$$\\

$$\Delta'':
\begin{array}{ccccccc}
    &  & \underline{e}_0'' & \rightarrow & \underline{e}_1'' &  &  \\
   & \quad\nearrow &  &  &  & \searrow & \\
   & \underline{e}_{w-1}'' &  &  &  & \quad \underline{e}_2'' &  \\
   & \uparrow &  &  &  & \quad\downarrow &  \\
   & \vdots &  &  &  & \quad\vdots &  \\
   & \quad \nwarrow&  &  &  & \swarrow &  \\
   &  & \underline{e}_{i+1}'' & \leftarrow & \underline{e}_i'' &  &  \\
\end{array}$$\\

\noindent Given a $C$-orbit of a simple regular module, the
corresponding tube consists of the indecomposable regular modules
whose regular composition factors belong to this orbit. We call
the tube corresponding to the orbit of a homogeneous module a
homogeneous tube.\\
In particular, in the sections that follow, we will need the following fact ([DR] Theorem 5.1):\\

\begin{lemma}
There exists a regular map $V:K^2\setminus\{(0,0)\}\rightarrow
Rep(Q,\underline{h})$ with the following properties:\\
\item{i)}\;$V(\varphi,\psi)$ is an indecomposable object in
$R_{(\varphi:\psi)}$ for each $(\varphi,\psi)\in
K^2\setminus\{(0,0)\}$ (it has to be a simple object in
$R_{(\varphi:\psi)}$
if $R_{(\varphi:\psi)}$ is a homogeneous tube).\\
\item{ii)}\;If $(\varphi:\psi)=(\gamma:\delta)$, then $V(\varphi,
\psi)$ and $V(\gamma,\delta)$ are isomorphic.\;$\diamondsuit$
\end{lemma}

\bigskip

\noindent Denote by $D_r$ the set of dimension vectors of all
regular representations of $Q$. Each element $\underline{d}\in
D_r$ can be written in the form
\begin{equation}
\underline{d} = p\underline{h}+\sum_{i\in I}
{p_i\underline{e}_i}+\sum_{i\in
I'}{p'_i\underline{e}'_i}+\sum_{i\in I''}{p''_i\underline{e}''_i}
\end{equation}
for nonnegative integers $p, p_i, p'_i, p''_i$ with at least one
coefficient in each family $(p_i\;|\;i\in I), (p'_i\; |\;i\in I'),
(p''_i \;|\;i\in I'')$ being zero. The decomposition in (1) is
called \textbf{canonical decomposition} of $\underline{d}\in D_r$
([Ri], Section 1). Since the only linear relations among the
dimension vectors $\underline{h},\underline{e}_i,
\underline{e}_i', \underline{e}_i''$ are the following
$$\underline{h} = \sum_{i\in I}{\underline{e}_i}=\sum_{i\in
I'}{\underline{e}_i'}=\sum_{i\in I''}{\underline{e}_i''},$$ we
have that such decomposition is unique. For simplicity, we set
\begin{equation}
\underline{d}'=\sum_I{p_i\underline{e}_i}+\sum_{I'}{p'_i\underline{e}'_i}+\sum_{I''}{p''_i\underline{e}''_i}.
\end{equation}

\noindent Finally, we recall the description of the generic decomposition of a regular dimension vector presented in [SW], Section 5.\\
We write
$\sum_I{p_i\underline{e}_i}=\sum_{j=1}^{max(p_i)}\sum_{k=1}^{N_j(p)}{\underline{v}_{j,k}}$,
where the dimension vectors are defined by reverse induction on
max($p_i$) as follows. Let $s=$ max($p_i$). Look at the set
$M_s=\{i\in\; I; p_i=s\}$ and decompose it $M_s=M_{s,1}\cup \cdots
\cup M_{s,t}$ into a sum of connected components (a subset of $I$
is called connected if it is an arc of the polygon $\Delta$). Take
$t=N_s(p)$ and $\underline{v}_{s,k}=\sum_{i\in\;
M_{s,k}}\underline{e}_i.$ Then repeat the procedure for
$\sum_I{\tilde{p}_i\tilde{\underline{e}}_i}=\sum_I{p_i\underline{e}_i}
- \sum_{k=1}^{N_s(p)}{\underline{v}_{s,k}}$, because we have
max$(\tilde{p_i})=s-1$. Similarly, we write
$\sum_I'{p'_i\underline{e}'_i}=\sum_{j=1}^{max(p'_i)}\sum_{k=1}^{N_j(p')}{\underline{v}'_{j,k}}$,
and
$\sum_I''{p''_i\underline{e}''_i}=\sum_{j=1}^{max(p''_i)}\sum_{k=1}^{N_j(p'')}{\underline{v}''_{j,k}}$.\\
\noindent In general, we have that the following decomposition:\\
\begin{equation}
 \underline{d} = \sum\limits_{i=1}^p\underline{h}
+\sum\limits_{j=1}^{max(p_i)}\sum\limits_{k=1}^{N_j(p)}{\underline{v}_{j,k}}
+
\sum\limits_{j=1}^{max(p'_i)}\sum\limits_{k=1}^{N_j(p')}{\underline{v}'_{j,k}}
+
\sum\limits_{j=1}^{max(p''_i)}\sum\limits_{k=1}^{N_j(p'')}{\underline{v}''_{j,k}}
\end{equation}
{\ }\\

\noindent is the generic decomposition of $\underline{d}$ ([SW], Proposition 44).\\

\section{Main results about algebras of semi-invariants}

\noindent Before giving the main theorem, we recall the following
notations that appear in [SW], Section 2. Let $Q$ be a Euclidean
quiver, $\underline{d} \in D_r$ with the canonical
decomposition as in (1) with $p\geq 1$.\\
We label the vertices $\underline{e}_i, \underline{e}_i',
\underline{e}_i''$ of the polygons $\Delta, \Delta', \Delta''$ in
Proposition 2.2, with the coefficients $p_i, p_i', p_i''$. Now, in
these new polygons, that we call $\Delta(\underline{d}),$
$\Delta'(\underline{d}),$ $\Delta''(\underline{d})$, we label the
edge from $p_k$ to $p_{k+1}$ with $\underline{e}_k$, for
$k=0,\ldots,u-2$, and the edge from $p_{u-1}$ to $p_0$ with
$\underline{e}_{u-1}$ (similarly for the other
polygons).\\
We say that the labeled arc $p_i --\cdots--p_j$ (in clockwise
orientation) of $\Delta(\underline{d})$ is $\textbf{admissible}$
if $p_i=p_j$ and $p_i<p_k$ for all its interior labels $p_k$. We
denote such arc by $[i,j]$. Similarly, we define admissible arcs
for the polygons $\Delta'(\underline{d})$ and
$\Delta''(\underline{d})$. Denote by $\mathcal{A}(\underline{d}),
\mathcal{A}'(\underline{d}), \mathcal{A}''(\underline{d})$ the
sets of all admissible labeled arcs in the polygons
$\Delta(\underline{d}),
\Delta'(\underline{d}), \Delta''(\underline{d})$.\\
On the other hand, we denote by $E_{i,j}$, for $i,j\in\{0,\ldots,
u-1\}$ with $j\neq i+1$ and $j\neq 0$ if $i=u-1$, the
in\-de\-com\-po\-sa\-ble regular module of dimension vector
$\underline{d}_{E_{i,j}}$ with the canonical decomposition as in
(1) with $p=0$, with socle $E_{i,i}$ and top $E_{j,j}$, where
$E_{k,k}$, or simply $E_k$ is the non homogeneous simple regular
module which has dimension vector $\underline{e}_k$, as in
Proposition 2.2. One gives analogous definitions for $E_{r,s}'$
and $E_{t,m}''$. Using the above notations, the dimension vector
of $E_{i,j}$ is equal to the sum of all $\underline{e}_k$ which
appear as edges in the polygon $\Delta(\underline{d})$ between the
vertices $p_j$ and $p_{i+1}$ (in clockwise orientation), where we
put $p_u:=p_0$
(similarly for $E_{r,s}'$ and $E_{t,m}''$).\\
Now we are ready to state the following theorem:

\begin{teo}
Let $Q$ be a Euclidean quiver, $\underline{d}\in D_r$. Let
$\underline{d}=p\underline{h}+\underline{d}'$ be the canonical
decomposition of \underline{d}, as in (1), with
$p\geq 1$.\\
Then the algebra $SI(Q,\underline{d})$ is generated as a
$K$-algebra by Schofield semi-invariants corresponding to the
indecomposable regular modules $E_{i,j}$, $E_{r,s}'$ and
$E_{t,m}''$ for each pair $(i,j), (r,s)$ and $(t,m)$ such that
$[j,i+1]
 \in \mathcal{A}(\underline{d})$,
 $[s,r+1]\in\mathcal{A}'(\underline{d})$ and $
[m,t+1]\in\mathcal{A}''(\underline{d})$ and by semi-invariants
$c_0, c_1, \ldots, c_p$ of weight $\partial$.
 The ideal of relations among generators is generated by the following relations:\\
$$
c_0=\prod{c^{E_{i,j}}},\; c_p=\prod{c^{E'_{r,s}}},\; c_0+\ldots+c_p = \prod{c^{E''_{t,m}}}\\
$$
where the products are over the pairs of indices $(i,j), (r,s),
(t,m)$ such that
$\sum{\underline{d}_{E_{i,j}}}=\sum{\underline{d}_{E'_{r,s}}}=\sum{\underline{d}_{E''_{t,m}}}=\underline{h}$,
respectively.
\end{teo}

\noindent \emph{Remark:} Our assumption about the dimension
vector, that is $\underline{d}\in D_r$ and $\underline{d}
\neq\underline{d}'$, is natural, because in the other cases
$SI(Q,\underline{d})$ is a polynomial algebra. This is due to the
following two facts ([SK] Section 4 Proposition 5; [Ri] Corollary
2.4, Theorem 3.2, Theorem 3.5)

\begin{lemma}\textbf{(Sato-Kimura)}
Suppose that $GL(\underline{\beta})$ has a dense orbit in
$Rep(Q,\underline{\beta})$. Let $S$ be the set of all $\sigma$
such that there exists an $f_\sigma \in
SI(Q,\underline{\beta})_\sigma$ which is nonzero and irreducible.
Then:
\item{a)} For every weight $\sigma \in S$, we have $dim SI(Q,\underline{\beta})_\sigma
\leq 1$.
\item{b)} All weights in $S$ are linearly independent over
$\mathbb{Q}$.
\item{c)} $SI(Q,\underline{\beta})$ is the polynomial ring generated by all
$f_\sigma, \sigma\in S$.\quad$\diamondsuit$
\end{lemma}

\begin{prop}
For a dimension vector $\underline{d}\in \mathbb{N}^{Q_0}$, the
variety $Rep\,(Q,\underline{d})$ has no open
$GL(\underline{d})$-orbit if and only $\underline{d}\in D_r$ and
for the canonical decomposition
$\underline{d}=p\underline{h}+\underline{d}'$ we have $p\geq 1$
(equivalently $\underline{d}
\neq\underline{d}'$).\quad$\diamondsuit$
\end{prop}
{\ }\\

\bigskip

\section{Proof of Theorem 3.1}

\noindent We divide the proof of Theorem 3.1 into two principal steps:\\

\noindent\emph{Step 1:\quad Description of generators}\\

     \noindent By Theorem 1.2, $SI(Q,\underline{d})$ is generated, as a $K$-algebra, by Schofield semi-invariants
      $\{c^V\}$ for all $V$ indecomposable representations of $Q$
    such that $\langle{\underline{d}_V,\underline{d}}\rangle=0$.
    Equivalently, by the definition of the generic decomposition of a
    dimension vector, we consider as generators $\{c^V\}$, for all $V$ indecomposable
    modules such that $\underline{d}_V$ is orthogonal to each summand of the generic decomposition
    of $\underline{d}$. Since vector $\underline{h}$ always appears among these
    summands,  we'll have that all $V$ belong to the set of regular
    indecomposable modules. Finally, by Proposition 1.1 e),
    we have  $$c^{V}\neq 0\; \Longleftrightarrow \; Hom_Q(V,W_i)= 0$$ for all $W_i$ summands of generic decomposition of
    $\underline{d}$, as in (3).\\
    \noindent Moreover,  by the definition of the generic decomposition of a
    regular dimension vector $\underline{d}$ as in (3), Proposition 1.1 c) and Lemma 3.2,
     we are reduced to consider as generators of  $SI(Q,\underline{d})$
    the following families:
an infinite family given by $\{c^{V(\varphi,\psi)}\}$, where
$(\varphi,\psi)\in\; K^2\setminus \{(0,0)\}$ and three finite
families given by $\{c^{E_{i,j}}\}$, $\{c^{E'_{r,s}}\}$ and
$\{c^{E''_{t,m}}\}$, where $E_{i,j}$, $E_{r,s}'$ and $E_{t,m}''$
are as in Section 3 and they satisfy two conditions
\begin{enumerate}

\item $\langle{\underline{d}_{E_{i,j}},\underline{d}}\rangle=0$, (similarly for $\underline{d}_{E_{r,s}'}$ and $\;\underline{d}_{E_{t,m}''}$);
\item $ Hom_Q(E_{i,j}, W_k)=0$ for all $W_k$ summands of the generic
decomposition of $\underline{d}$, (similarly for $E_{r,s}'$ and
$E_{t,m}''$).
\end{enumerate}
{\ }\\

Now, remembering the notation in Section 3, we want to prove that
the above conditions 1 and 2 are equivalent to the condition of
admissibility of the arc $[j,i+1]$ in the polygon
$\Delta(\underline{d})$
($[s,r+1]\in \Delta'(\underline{d})$ and $[m,t+1]\in \Delta''(\underline{d})$).\\
We start proving the following lemma
\begin{lemma}
In the above notation, the condition 1, that is
$\langle{\underline{d}_{E_{i,j}},\underline{d}}\rangle=0$,
 is equivalent to consider the arc
$[j,i+1]$ in the polygon $\Delta(\underline{d})$ with extremes
$p_j$ and $p_{i+1}$ that are the same number.
\end{lemma}

\emph{Proof:} As in Section 3, we write $\underline{d}_{E_{i,j}}$
as the sum of all $\underline{e}_k$ which label the edges in the
arc $[j, i+1]$. Recall the following result ([DR], Lemma 3.3):

\begin{lemma}Let $\underline{e}_i, \underline{e}_j\; (\underline{e}_k', \underline{e}_r'$ and $\underline{e}_m'',
\underline{e}_n''$ respectively) two dimension vectors as in
Proposition 2.2, then

    $${\langle{\underline{e}_i,\underline{e}_j}\rangle}=\bigg\{%
\begin{array}{ll}
    \;1 & \hbox{if i\,=\,j;} \\
    -1 & \hbox{if i\,=}\,j-1; \\
    \;0 & \hbox{otherwise}; \\
\end{array}%
$$ where $\underline{e}_{-1}:=\underline{e}_{u-1}$ (we have analogous result for the
other pairs). Otherwise, the value of Euler form between dimension
vectors belonging to different polygons of Proposition 2.2, is
zero.\qquad$\diamondsuit$
\end{lemma}

\noindent We put $\underline{d}_{E_{i,j}}:=
\sum\limits_{k=0}^{u-1}{\delta_k\underline{e}_k}$, where
$\delta_k$ is equal to 1, for each $k$ such that
$\underline{e}_k$ is an edge in $[j,i+1]$, and zero otherwise. Using Lemma 4.2, we have the following equalities:\\

$$\langle{\underline{d}_{E_{i,j}},\underline{d}}\rangle=
\langle{\sum_{k=0}^{u-1}{\delta_k\underline{e}_k},
    p\underline{h}+\sum_{i\in I}{p_i\underline{e}_i}+\sum_{j\in I'}{p_j'\underline{e}_j'}+\sum_{k\in I''}{p_k''\underline{e}_k''}}\rangle
    =$$
$$=\langle{\sum_{k=0}^{u-1}{\delta_ke_k},
    \sum_i{p_i\underline{e}_i}+\sum_j{p_j'\underline{e}_j'}+\sum_k{p_k''\underline{e}_k''}}\rangle =$$

    $$=\sum_k{\delta_kp_k\langle{\underline{e}_k,\underline{e}_k}\rangle + \sum_{k\neq
    l}{\langle{\delta_k\underline{e}_k, p_l\underline{e}_l}\rangle}=\sum_k{\delta_kp_k}+\sum_{k\neq l}{\delta_kp_l\langle{\underline{e}_k,
    \underline{e}_l}\rangle}}=$$
    $$=\sum_k{\delta_kp_k}+\sum_{l=0}^{u-1}{\delta_{l-1}p_l\langle{\underline{e}_{l-1},
    \underline{e}_l}\rangle}=
    \sum_k{\delta_kp_k}-\sum_{l=0}^{u-1}{\delta_{l-1}p_l}=$$
    $$=p_j-p_{i+1}.$$\\

\noindent In conclusion, we obtain
\begin{equation}
\langle{\underline{d}_{E_{i,j}},\underline{d}}\rangle
=0\;  \Longleftrightarrow\; p_j=p_{i+1}\end{equation}\\
\noindent as requested.\qquad$\diamondsuit$\\

\noindent Finally, we prove that the conditions 1 and 2 together
are equivalent to the condition of admissibility of the arc
$[j,i+1]\in \mathcal{A}(\underline{d})$. From now on, we'll use
the short notation, $n\in [i, j]$, to say that $n$ is such that
$\underline{e}_n$ is an edge in the oriented arc $[i,j]$ of
$\Delta(\underline{d})$.\\
By the definition of admissibility and Lemma 4.1, we are reduced
to
prove the following statement:\\

\begin{lemma}Let $p_j,p_{i+1}\in \Delta(\underline{d})$ such that $p_j=p_{i+1}$. Then there exists
$l\in [j+1,i+1]$ such that $p_l<p_j$ if and only if there exists
$\tilde{W}$ summand of the generic decomposition of \underline{d}
such that $Hom_Q(E_{i,j},\tilde{W})\neq 0$.
\end{lemma}
{\ }\\

\noindent Before starting the proof of Lemma 4.3, we recall a
simple fact that follows by the structure of Auslander-Reiten
quiver of extended Dynkin quivers ([SS], Chapter XIII.2 Theorem
2.1):
\begin{lemma}
Let $E_{k,r}$ and $E_{m,n}$ be two regular indecomposable non
homogeneous modules of $Q$ as in Section 3. Then  $Hom_Q(E_{k,r},
E_{m,n})\neq 0$ if and only if $E_{m,n}$ with $m\in [r,k+1]$ and
$n\in [m+2,r+1]$.\qquad$\diamondsuit$
\end{lemma}
{\ }\\

\noindent\emph{Proof of Lemma 4.3:} Let $p_l$ be the minimal
number in the arc $[j,i+1]\in \Delta(\underline{d})$
that is nearest to $p_j$. Graphically we have\\
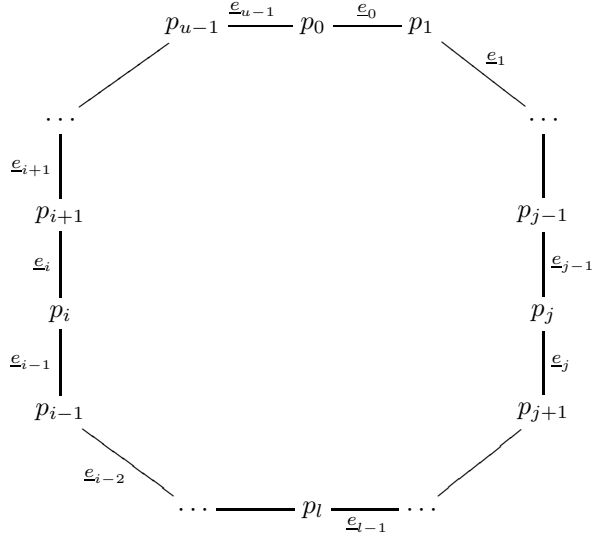
$$
\begin{xy}
\xymatrix{ &p_{u-1}\ar@{-}[r]^{\underline{e}_{u-1}}&p_0\ar@{-}[r]^{\underline{e}_0}&p_1\ar@{-}[dr]^{\underline{e}_1}\\
\cdots\ar@{-}[ur]\ar@{-}_{\underline{e}_{i+1}}[d]& & & &\cdots\ar@{-}[d]\\
p_{i+1}\ar@{-}_{\underline{e}_i}[d]& & & &p_{j-1}\ar@{-}[d]^{\underline{e}_{j-1}}\\
p_{i}\ar@{-}_{\underline{e}_{i-1}}[d]& & & &p_j\\
 p_{i-1}\ar@{-}_{\underline{e}_{i-2}}[dr]& & & &p_{j+1}\ar@{-}[u]_{\underline{e}_j}\\
 &\cdots\ar@{-}[r]&p_l\ar@{-}[r]_{\underline{e}_{l-1}}&\cdots\ar@{-}[ur]\\
 }
\end{xy}\\
$$

\bigskip

\noindent By definition of the generic decomposition of
$\underline{d}$ (see (3)), there exists a summand of the generic
decomposition of $\underline{d}$, $\tilde{W}$, which is isomorphic
to $E_{l-1,n}$, with $n\in [l+1,
j+1]$. By Lemma 4.4, we have $ Hom_Q(E_{i,j}, \tilde{W})\neq0$, as requested.\\
Viceversa, if we have a summand of the generic decomposition of
$\underline{d}$, $\tilde{W}$, such that
$Hom_Q(E_{i,j},\tilde{W})\neq 0$, then, by Lemma 4.4, it is
isomorphic to $E_{m,n}$, with $m\in [j,i+1]$ and $n\in [m+2,
j+1]$. Now, by the hypothesis that $p_{i+1}=p_j$, no summand of
the generic decomposition of $\underline{d}$ can be isomorphic to
$E_{i,r}$ with $r\in [i+2, j+1]$. In all other cases, by the
construction of the generic decomposition as in (3), it follows
that $p_{m+1}<p_j$ and $m+1\in
[j+1,i+1]$, as requested.\quad$\diamondsuit$\\

\bigskip

\noindent By Lemma 3.2, Proposition 3.3 and the reciprocity
formula ([DW] Corollary 1), we have that
$$dim \;SI(Q,
\underline{d})_{\langle{\underline{d}_{E_{i,j}},-\rangle}}=dim
\;SI(Q,
\underline{d})_{\langle{\underline{d}_{E'_{r,s}},-\rangle}}=dim
\;SI(Q,
\underline{d})_{\langle{\underline{d}_{E''_{t,m}},-\rangle}}=1$$
for each pair $(i,j), (r,s)$ and $(t,m)$ such that $[j,i+1]
 \in \mathcal{A}(\underline{d})$,
 $[s,r+1]\in\mathcal{A}'(\underline{d})$ and $
[m,t+1]\in\mathcal{A}''(\underline{d})$), so each semi-invariant
$c^{E_{i,j}}$, $c^{E'_{r,s}}$ and $c^{E''_{t,m}}$, where $(i,j)$,
$(r,s)$ and $(t,m)$ are as above, spans the corresponding weight
space. All the other generators, $c^{V(\varphi,\psi)}$ with
$(\varphi, \psi)\in\; K^2\setminus \{(0,0)\}$ are in $SI(Q,
\underline{d})_{\partial}$. By the next result, it follows that
$dim \;SI(Q,
\underline{d})_{\partial}=p+1$.\\

\begin{prop}
If \;$\underline{d}=p\underline{h}+\underline{d}'$ is the
canonical decomposition of $\underline{d}$ as in (1), then $dim \;
SI(Q,\underline{d})_{m\partial}={p+m \choose m}$.
\end{prop}
{\ }\\

\noindent\emph{Proof}: Applying the reciprocity formula, we have
$$dim \;SI(Q,\underline{d})_{m\partial}=dim
\;SI(Q,m\underline{h})_{-\langle{-,\underline{d}}\rangle}.$$

Let $c_W$ be an element of
$SI(Q,m\underline{h})_{-\langle{-,\underline{d}}\rangle}$, where
$W$ is a representation of $Q$ of dimension vector
$\underline{d}$. Using the generic decomposition of
$\underline{d}$ as in (3) and Proposition 1.1 d), $c_W=c_{W'}
c_{W''}$ where $W'$ and $W''$ are representations of dimension
vectors $p\underline{h}$ and $\underline{d'}$, respectively.
\noindent By Lemma 3.2 and Proposition 3.3, it follows that
$dim\;SI(Q,m\underline{h})_{-\langle{-,\underline{d'}}\rangle}=1$,
so we have an isomorphism between
$SI(Q,m\underline{h})_{-\langle{-,\underline{d}}\rangle}$ and $
SI(Q,m\underline{h})_{-\langle{-,p\underline{h}}\rangle}$ given by
multiplication by a non-zero semi-invariant of weight
$-\langle{-,\underline{d'}}\rangle$.\\
We need to show $$dim\;
SI(Q,m\underline{h})_{-\langle{-,p\underline{h}}\rangle}= dim\; SI(Q,p\underline{h})_{\langle{m\underline{h},-}\rangle}={p+m \choose m}.$$\\
Applying Theorem 1.3, we are reduced to count all
subrepresentations of dimension vector $m\underline{h}$ in a
general representation of dimension vector $(p+m)\underline{h}$.
Since a general representation of dimension vector
$(p+m)\underline{h}$ is a direct sum of $(p+m)$ different
indecomposables of dimension vector $\underline{h}$, then it has
${p+m \choose m}$ subrepresentations
of dimension vector $m\underline{h}$.\qquad$\diamondsuit$\\

\bigskip

\noindent\emph{Step 2:\quad Description of relations}\\

\noindent We conclude the proof of Theorem 3.1 giving a
description of the relations among the generators of
$SI(Q,\underline{d})$ presented in Step 1. In particular, in the
proof of the following proposition, we show that the set of
semi-invariants $\{c_0,\ldots,c_p\}$ that appears in the statement
of Theorem 3.1, is a basis of $SI(Q,\underline{d})_\partial.$ We
recall the following fact of linear algebra that we use in the
proof of Proposition 4.7:

\begin{lemma}
A \textbf{Vandermonde} matrix is a $m\times n$ $K$-matrix,
$M=(m_{i,j})$ where

$$m_{i,j}=(\alpha_{j})^{i-1},$$

\noindent for all indices $i$ and $j$. If $m=n$, the determinant
of a square Vandermonde matrix can be expressed as:

$$
det(M)=\prod_{1\leq i< j\leq n}(\alpha_j-\alpha_i).
$$

\noindent Thus, a square Vandermonde matrix is invertible if and
only if the $\alpha_i$ are di\-stinct.
\end{lemma}

\begin{prop}
\item{(a)} We have the following relations among the generators
of the algebra $SI(Q,\underline{d})$
\begin{equation}
c_0=\prod{c^{E_{i,j}}},\qquad c_p=\prod{c^{E'_{r,s}}},\qquad
c_0+\ldots+c_p = \prod{c^{E_{t,m}''}},
\end{equation}
where the products are over the pairs of indices $(i,j), (r,s),
(t,m)$ such that
$\sum{\underline{d}_{E_{i,j}}}=\sum{\underline{d}_{E'_{r,s}}}=\sum{\underline{d}_{E''_{t,m}}}=\underline{h}$,
respectively.
\item{(b)} The relations in (5) are enough to generate the ideal of relations.
\end{prop}
{\ }\\

\noindent \emph{Proof:}
 \emph{(a)}\;By Step 1, we have that Schofield semi-invariants $c^{V(\varphi,\psi)}$, with $(\varphi,\psi)\in K^2\setminus \{(0,0)\}$, span
$SI(Q,\underline{d})_{\partial}$. Let
$\underline{d}=p\underline{h}+\underline{d}'$ be the canonical
decomposition of $\underline{d}$ as in (1). Using the same
arguments of the proof of Proposition 4.5, we have that
$c^{V(\varphi,\psi)}\in SI(Q,\underline{d})_{\partial}$ is equal
to the product of the semi-invariant
$\tilde{c}^{V(\varphi,\psi)}:=c^{V(\varphi,\psi)}\in\;SI(Q,p\underline{h})_{\partial}$
and a non-zero semi-invariant
$f'\in SI(Q,\underline{d'})_{\partial}$.\\
We want to focus on the first factor $\tilde{c}^{V(\varphi,\psi)}$.\\
For a fixed $V(\varphi,\psi)$ $\in Rep(Q, \underline{h})$,
$\tilde{c}^{V(\varphi,\psi)}$ is the restriction of the function
$\tilde{c}:=c\in K[Rep(Q, \underline{h})]\otimes
K[Rep(Q,p\underline{h})]$ to
$\{V(\varphi,\psi)\}\times Rep (Q,p\underline{h})$.\\
On the other hand, by definition of the generic decomposition in
(3), a general representation of dimension vector $p\underline{h}$
is a direct sum of $p$ pairwise non isomorphic representations
$V(\gamma_i,\delta_i)$,
$i=1,\ldots,p$ of dimension vectors $\underline{h}$.\\
Thus, if we define the following map

\begin{align}
K^2\times (K^2)^p&\stackrel{\pi}\rightarrow
Rep(Q,\underline{h})\times
Rep(Q,p\underline{h})\\
\nonumber(\varphi, \psi,
\gamma_1,\delta_1,\ldots,\gamma_p,\delta_p)&\mapsto
V(\varphi,\psi)\oplus \bigoplus_{i=1}^p{V(\gamma_i,\delta_i)},
\end{align}

\noindent we want to study the image in $K[\varphi, \psi]\otimes
K[\gamma_i, \delta_i\mid i=1,\ldots,p]$ of the space
$SI(Q,p\underline{h})_{\partial}=Span_K\{\tilde{c}^{V(\varphi,\psi)}\mid
(\varphi,\psi)\in\;K^2\setminus \{(0,0)\}\}$ under the map $\pi^*$ induced by (6).\\

\noindent First of all, if $W$ is a general representation of
dimension vector $p\underline{h}$, applying Proposition 1.1 d), we
have that:

$$\tilde{c}^{V(\varphi,\psi)}(W)=\tilde{c}_W(V(\varphi,\psi))=\prod_{i=1}^p{\hat{c}_{V(\gamma_i,\delta_i)}(V(\varphi,\psi))}=
\prod_{i=1}^p{\hat{c}^{V(\varphi,\psi)}(V(\gamma_i,\delta_i))},$$

\noindent where $\hat{c}:=c\in K[Rep(Q, \underline{h})]\otimes
K[Rep(Q,\underline{h})]$. Thus, we are reduced to study the case $p=1$. We need to prove the following fact:\\

\begin{lemma}
Let $\hat{c}^{V(\varphi,\psi)}\in\;
SI(Q,\underline{h})_{\partial}$. Then:\\

\item{i)}\;For each $(\varphi,\psi)$ and $(\gamma, \delta)\in K^2\setminus \{(0,0)\}$, we have that

$$\hat{c}^{V(\varphi,\psi)}(V(\gamma,\delta))=0\Longleftrightarrow(\varphi:\psi)=(\gamma:\delta).$$

\item{ii)}\;The image of $\hat{c}$ in $K[\varphi,
\psi]\otimes K[\gamma, \delta]$ under the map $\pi^*$ induced by
(6) is a polynomial of the form $\varphi\delta-\psi\gamma$.
\end{lemma}

\noindent\emph{Proof: i)} follows by Proposition 1.1 e) and the
following consequence of [DR], Theorem 3.5:
$$Hom_Q(V(\varphi,\psi),V(\gamma,\delta))\neq 0\Longleftrightarrow
(\varphi:\psi)=(\gamma:\delta).$$

\noindent \emph{ii)} First of all, by \emph{i)} follows that the
image of $\hat{c}$ in $K[\varphi, \psi]\otimes K[\gamma, \delta]$
is a polynomial of the form
$$(\varphi\delta-\psi\gamma)^m=\sum_{j=0}^m {{m \choose
j}(\varphi\delta)^{m-j}(-\psi\gamma)^j},$$ $m\geq 1$, where the
pairs of variables $(\varphi, \psi)$ and $(\gamma,\delta)$ are two
generators of a homogeneous coordinate ring of $\mathbb{P}_1(K)$,
so they are linearly independent, respectively. We want to prove that $m=1$.\\
Suppose $m>1$. Then, the dimension of the image in $K[\varphi,
\psi]\otimes K[\gamma, \delta]$ of the vector space
$SI(Q,\underline{h})_{\partial}=Span_K\{\hat{c}^{V(\varphi,\psi)}\mid
(\varphi,\psi)\in\;K^2\setminus \{(0,0)\}\}$ would be bigger than
2, a contradiction since by Proposition 4.5 we have that
$dim_K\;SI(Q,\underline{h})_{\partial}=2$. Indeed, for example,
the image of $Span_K\{\hat{c}^{V(\varphi,\psi)}\mid
(\varphi,\psi)\in\;K^2\setminus \{(0,0)\}\}$ contains the vector
space $$P:=Span_K\{\delta^m,\gamma^m,\sum_{j=0}^m {{m \choose
j}(\delta)^{m-j}(-\gamma)^j}\},$$
\noindent which has dimension equal to 3.\qquad${\diamondsuit}$\\

\bigskip

\noindent Now, if $p>1$, by Lemma 4.7 \emph{ii)}, we have that the
image of $\tilde{c}$ is equal to the polynomial
$$f:=\prod_{i=1}^p(\varphi\delta_i-\psi\gamma_i).$$\\
Without losing generality, we put $\delta_i=1$ for all $i$
(similarly if we put $\gamma_i=1$). Then, we have that
$$
f=a_0(\gamma_1,\ldots,\gamma_p)\varphi^p+a_1(\gamma_1,\ldots,\gamma_p)\varphi^{p-1}\psi+\ldots+a_p(\gamma_1,\ldots,\gamma_p)\psi^p,
$$
where $a_i(\gamma_1,\ldots,\gamma_p)$ are the elementary symmetric
functions in $\gamma_1,\ldots,\gamma_p$.\\
In particular, the image of $SI(Q,p\underline{h})_{\partial}$ is
equal to $Span_K\{f_{(\overline{\varphi},\overline{\psi})},
(\overline{\varphi},\overline{\psi})\in K^2\setminus \{(0,0)\}\}$,
where $f_{(\overline{\varphi},\overline{\psi})}\in K[\gamma_i,
\delta_i\mid i=1,\ldots,p]$ is the valuation of $f$ in $(\varphi,\psi)=(\overline{\varphi},\overline{\psi})$.\\
Moreover, we have that
$Span_K\{f_{(\overline{\varphi},\overline{\psi})},
(\overline{\varphi},\overline{\psi})\in K^2\setminus \{(0,0)\}\}=Span_K\{a_0,\ldots,a_p\}$.\\
\noindent Indeed, we have that
$Span_K\{f_{(\overline{\varphi},\overline{\psi})},
(\overline{\varphi},\overline{\psi})\in K^2\setminus
\{(0,0)\}\}\subseteq Span_K\{a_0,\ldots,a_p\}$, so we need to
prove the other inclusion. Without losing generality we put
$\overline{\psi}=1$. It's sufficient to show that there exist
$\underline{k}'=(k'_0,\ldots,k'_p)\in K^{p+1}$ such that the
following equality holds:

$$
\sum_{i=0}^p{k_ia_i}=\sum_{i=0}^p{k'_i\overline{\varphi}_i^pa_i},
$$

\noindent where $k_i, \overline{\varphi}_i, k_i'\in K$ and
$\overline{\varphi}_i$ are pairwise distinct. It is equivalent to
solve the linear system $\underline{k}=A\underline{k}'$, where
$\underline{k}=(k_0,\ldots,k_p)\in K^{p+1}$, and $A$ is a square
Vandermonde matrix (see Lemma 4.6). Since $\overline{\varphi}_i$
are pairwise distinct, the matrix $A$ is invertible and the system
$\underline{k}=A\underline{k}'$
is compatible.\\

\noindent In conclusion, remembering that the functions $a_i$ for
$i=0,\ldots,p$ are li\-near\-ly independent, if we define
$c_i':=(\pi^*)^{-1}(a_i)$ for $i=0,\ldots,p$, we have that
$\{c_0',\ldots,c_p'\}$ is a basis of
$SI(Q,p\underline{h})_{\partial}$ as well as $\{c_i:=c_i'f',
i=0,\ldots,p\}$ is a basis of
$SI(Q,\underline{d})_{\partial}$.\\
In particular, we note that

$$c_0:=(\pi^*)^{-1}(a_0)f'=(\pi^*)^{-1}(f_{(1,0)})f'=c^{V(1,0)},$$

\noindent(similarly for $c_p$ and $c_0+\ldots+c_p$ that are the
Schofield semi-invariants corresponding to $V(0,1)$ and $V(1,1)$,
respectively). It's known (see [DR] Section 5 and 6) that the
modules $V(1,0),$ $V(0,1),$ $V(1,1)$ are consecutive extensions of
$E_{i,j}, E'_{r,s}, E''_{t,m}$, respectively, such that
$\sum{\underline{d}_{E_{i,j}}}=\underline{h}$,
$\sum{\underline{d}_{E'_{r,s}}}=\underline{h}$ and
$\sum{\underline{d}_{E''_{t,m}}}=\underline{h}$, respectively. Applying the Proposition 1.1 c), we obtain the relations in (5).\\

\noindent \emph{(b)} By Step 1 and point \emph{(a)}, we have that
the algebra $SI(Q,\underline{d})$ is of the form:

$$\frac{K[c_0,c_1,\ldots,c_p,c^{E_{i,j}}, c^{E'_{r,s}}, c^{E_{t,m}''}]}{I},$$

\noindent where $I$ is an ideal which contains the ideal $J$
generated by the relations in (5). To prove that the two ideals
$J$ and $I$ are the same ideal, we need to show that the
epimorphism $\mu$
$$
T:=\frac{K[c_0,c_1,\ldots,c_p,c^{E_{i,j}}, c^{E'_{r,s}},
c^{E_{t,m}''}]}{J}\buildrel{\mu}\over\longrightarrow
SI(Q,\underline{d})
$$
is an isomorphism.\\
We are going to show that the dimensions of weight spaces in both
rings are equal. If we look at $T$ and $SI(Q,\underline{d})$ as
$\bigoplus_{\sigma} {T_\sigma}$ and $
\bigoplus_{\sigma}{SI(Q,\underline{d})_\sigma}$ respectively, we
recognize an epimorphism between the corresponding graded
components
$$
T_\sigma\longrightarrow SI(Q,\underline{d})_\sigma.
$$
Then, we have
\begin{equation}
dim T_\sigma\geq dim SI(Q,\underline{d})_\sigma.
\end{equation}\\

\noindent  Now, using the same arguments of the proof of
Proposition 4.5, we have that
$$ SI(Q,\underline{d})_{\langle{\underline{\alpha},-}\rangle}\cong
SI(Q,\underline{d})_{\langle{m\underline{h},-}\rangle},$$\\
where $\underline{\alpha}$ is a regular dimension vector with the
canonical decomposition (see (1))
$\underline{\alpha}=m\underline{h}+\underline{\alpha}'$ and
$\underline{d}$ has the canonical decomposition
$\underline{d}=p\underline{h}+\underline{d}'$.\\
Then, applying Proposition 4.5 we have that

\begin{equation}
dim\;SI(Q,\underline{d})_{\langle{\underline{\alpha},-}\rangle}={p+m
\choose m},
\end{equation}

\noindent and, by (7) and (8), we obtain the inequality:
\begin{equation}
dim \;T_{\langle{\underline{\alpha},-}\rangle}\geq dim
\;SI(Q,\underline{d})_{\langle{\underline{\alpha},-}\rangle}={p+m
\choose m}.
\end{equation}\\

\noindent Thus, if we prove that

\begin{equation} dim \;T_{\langle{\underline{\alpha},-}\rangle}\leq {p+m \choose m},
\end{equation}\\

\noindent by (9) and (10) follows that $\mu$ is an isomorphism.\\

\noindent First of all, we observe that $f\in
T_{\langle{\underline{\alpha},-}\rangle}$ if and only if $f$ is a
polynomial equal to the following one

$$
\sum_{(\underline{j},n_{i,j},n_{r,s}',n_{t,m}'')}
{k_{\underline{j},n_{i,j},n_{r,s}',n_{t,m}''}} c_0^{j_0}\cdots
c_p^{j_p} \prod_{(i,j)}(c^{E_{i,
j}})^{n_{i,j}}\prod_{(r,s)}(c^{E_{r,s}'})^{n'_{r,s}}\prod_{(t,m)}(c^{E_{t,m}''})^{n''_{t,m}},$$
where $k_{\underline{j},n_{i,j},n_{r,s}',n_{t,m}''}\in K$ and the
sum is over all $(\underline{j},n_{i,j},n_{r,s}',n_{t,m}'')$ such
that
$$
\underline{\alpha}\buildrel{\star}\over=(\sum_{s=0}^p{j_s})\underline{h}+\sum_{(i,j)}{n_{i,j}
\underline{d}_{E_{i,j}}}+\sum_{(r,s)}{n_{r,s}'\underline{d}_{
E_{r,s}'}}+\sum_{(t,m)}{n_{t,m}''\underline{d}_{E_{t,m}''}}.
$$\\

\noindent Using the relations (5), we transform the polynomial
$f$ into the following one:\\

$$\sum_{(\underline{i},m_{i,j},m_{r,s}',m_{t,m}'')}
k_{\underline{i},m_{i,j},m_{r,s}',m_{t,m}''} c_0^{i_0}\cdots
c_p^{i_p}
\prod_{(i,j)}(c^{E_{i,j}})^{m_{i,j}}\prod_{(r,s)}(c^{E_{r,s}'})^{m'_{r,s}}\prod_{(t,m)}(c^{E_{t,m}''})^{m''_{t,m}},$$

\noindent with $i_0+\ldots+i_p = m'$, that is equivalent to change
the decomposition $(\star )$ into the following

$$
\underline{\alpha}\buildrel{\star\star}\over=\sum_{l=1}^{m'}\underline{h}+\sum_{(i,j)}{m_{i,j}
\underline{d}_{E_{i,j}}}+\sum_{(r,s)}{m_{r,s}'\underline{d}_{
E_{r,s}'}}+\sum_{(t,m)}{m_{t,m}''\underline{d}_{E_{t,m}''}}.
$$\\

By the uniqueness of the canonical decomposition, $m'$ is equal to
$m$, the multiplicity of $h$ in the canonical decomposition of
$\underline{\alpha}$. Moreover, we know that the only linear
relations among the dimension vectors
$\underline{h},\underline{e}_i, \underline{e}_i',
\underline{e}_i''$ in the canonical decomposition (1) are the
following: $\underline{h} = \sum_{i\in
I}{\underline{e}_i}=\sum_{i\in I'}{\underline{e}_i'}=\sum_{i\in
I''}{\underline{e}_i''},$ then we have that the weights of the
generators of $T$
 are linearly
independent except the relations
$$\underline{h} =
\sum_{(i,j)}{\underline{d}_{E_{i,j}}}=\sum_{(r,s)}{\underline{d}_{
E_{r,s}'}}=\sum_{(t,m)}{\underline{d}_{E_{t,m}''}}.$$ In
conclusion, the decomposition $(\star\star)$ is unique and $f$ is
equal to:
$$\sum_{\underline{i}}
k_{\underline{i}} c_0^{i_0}\cdots c_p^{i_p}\big(
\prod_{(i,j)}(c^{E_{i,j}})^{m_{i,j}}\prod_{(r,s)}(c^{E_{r,s}'})^{m'_{r,s}}\prod_{(t,m)}(c^{E_{t,m}''})^{m''_{t,m}}\big).$$

\noindent Thus, we have that
\begin{equation} dim \;T_{\langle{\alpha,-}\rangle}\leq {p+m \choose m},
\end{equation}\\

\noindent as requested.\qquad$\diamondsuit$\\

\end{document}